\documentclass[12pt]{amsart}
\usepackage{times}
\usepackage{amssymb}
\usepackage{amscd}
\def\Aff{\mathop{\mathrm {Aff}}\nolimits}
\def\aff{\mathop{\mathrm {aff}}\nolimits}

\def\Lie{\mathop{\mathrm {Lie}}\nolimits}
\def\Ln{\mathop{\mathrm {Ln}}\nolimits}
\def\tr{\mathop{\mathrm {tr}}\nolimits}

\newtheorem{theorem}{Theorem}[section]
\newtheorem{proposition}[theorem]{Proposition}
\newtheorem{lemma}[theorem]{Lemma}
\newtheorem{corollary}[theorem]{Corollary}
\newtheorem{remark}[theorem]{Remark}

\begin{document}
\title[Quantum Co-Adjoint Orbits...]{Quantum Co-Adjoint Orbits of the Group of Affine Transformations of the Complex Straight Line}
\author{Do Ngoc Diep and Nguyen Viet Hai}
\date{Version of August 10, 1999}
\address{Institute of Mathematics, National Center for Natural Sciences and Technology, P. O. Box 631, Bo Ho, 10.000, Hanoi, Vietnam}
\email{dndiep@hn.vnn.vn}
\thanks{This work was supported in part by the National  Foundation for Fundamental Science Research of Vietnam and the Alexander von Humboldt Foundation, Germany}
\maketitle
\begin{abstract}
We construct start-products on the co-adjoint orbit of the Lie group
$\Aff({\bf C})$ of affine transformations of the complex straight line and apply them to obtain the irreducible unitary representations of this group. These results show effectiveness of the Fedosov quantization even for groups which are neither nilpotent nor exponential. Together with the result for the group  $\Aff({\bf R})$ [see DH], we have thus a description of quantum $\overline{MD}$ co-adjoint orbits.  
\end{abstract}
\section{Introduction}

The notion of $\star$-products was a few years ago introduced  and played a fundamental role in the basic problem of quantization, see e.g. references [AC1,AC2,F,G\dots], as a new approach to quantization on arbitrary symplectic manifolds. In [DH] we have constructed star-products on upper half-plane, obtained operator $\hat{\ell}_z,  z \in \aff({\bf R}) = \Lie\Aff({\bf R})$ and we have proved that the representation $\exp(\hat{ \ell}_z) = \exp(\alpha\frac{\partial}{\partial s} + i{\beta}e^s)$ of group $\Aff_0({\bf R})$ coincides with representation $T_{\Omega_\pm}$  obtained from the orbit method or Mackey small subgroup method. One of the advantage of this group, with which the computation is rather accessible is the fact that its connected component $\Aff_0(\mathbf R)$ is exponential. We could use therefore the canonical coordinates for  Kirillov form on the orbits.

It is natural to consider the same problem for the group $\Aff({\bf C})$. We can expect that the calculations and final expressions could be similar to the corresponding in real line case, but this group $\Aff({\bf C})$ is no more exponential, i.e exponential map $$\exp : \aff({\bf C}) \to \widetilde{\Aff}({\bf C})$$  is no longer a global diffeomorphism
and the general theory of D. Arnal and J. Cortet \cite{arnalcortet1}, \cite{arnalcortet2} and others could not be directly applicable. 
We overcame these difficulties by a way rather different which could indicate new ideas for more general non-exponential groups: To overcome the main difficulty in applying the deformation quantization to this group, we replace the global diffeomorphism in Arnal-Cortet's setting by a local diffeomorphism. With this replacement, we need to pay attention on complexity of the  the symplectic Kirillov form in new coordinates. We then computed the inverse image of the Kirillov form on appropriate local charts. The question raised here is how to choose a good local chart in order to have as possible simple calculation. The calculations we proposed are realized by using complex analysis on very simple complex domain. 

Our main result is the fact that by an exact computation we had found out an explicit star-product formula (Proposition 3.5) on the local chart. This means that the functional algebras on co-adjoint orbits admit a suitable deformation, or in other words, we obtained quantum co-adjoint orbits of this group as exact models of new quantum objects, say ``quantum punctured complex planes" $(\mathbf C^2 \setminus L)_q$. Then, by using the Fedosov deformation quantization, it is not hard to obtain the full list of irreducible unitary representations (Theorem 4.2) of the group $\Aff({\bf C})$, although the computation in this case, like by using the Mackey small subgroup method or modern orbit method, is rather delicate. The infinitesimal generators of those exact model of  infinite dimensional irreducible unitary representations, nevertheless,  are given by rather simple formulae.

We introduce some preliminary result in \S2. The operators $\hat {\ell}_A$ which define the representation of the Lie algebra $\aff({\bf C})$ are found in \S3. In particular, we obtain the unitary representations of Lie group $\widetilde{\Aff}$({\bf C}) in Theorem 4.3 \S4.

\section{Preliminary Results}
Recall that the Lie algebra ${\mathfrak g} = \aff({\bf C})$ of affine transformations of the complex straight line is described as follows, see [D].

It is well-known that the group $\Aff({\bf C})$ is a four (real) dimensional Lie group which is isomorphism to the group of matrices:
$$\Aff({\bf C}) \cong  \left\{\left(\begin{array}{cc} a & b \\ 0 & 1 \end{array}\right) \vert a,b \in {\bf C}, a \ne 0 \right\}$$

The most easy method is to consider $X$,$Y$ as complex generators,
$X=X_1+iX_2$ and $Y=Y_1+iY_2$. Then from the relation $[X,Y]=Y$, we get$ [X_1,Y_1]-[X_2,Y_2]+i([X_1Y_2]+[X_2,Y_1]) = Y_1+iY_2$. 
This mean that the Lie algebra $\aff({\bf C})$ is a real 4-dimensional Lie algebra, having 4 generators with the only nonzero Lie brackets: $[X_1,Y_1] - [X_2,Y_2]=Y_1$; $[X_2,Y_1] + [X_1,Y_2] = Y_2$ and we can choose another basic noted again by the same letters to have more clear Lie brackets of this Lie algebra:
$$[X_1,Y_1] = Y_1; [X_1,Y_2] = Y_2; [X_2,Y_1] = Y_2; [X_2,Y_2] = -Y_1$$

\begin{remark}{\rm
The exponential map $$\exp: {\bf C}  \longrightarrow  {\bf  C}^{*} := {\bf C} \backslash \{0\}$$  giving by $z \mapsto e^z$ is just the covering map and therefore the universal covering of ${\bf} C^*$ is $\widetilde {\bf C}^* \cong {\bf C}$. As a consequence one deduces that $$\widetilde {\Aff}({\bf C}) \cong {\bf C} \ltimes{\bf C} \cong  \{(z,w) \vert z,w \in {\bf C} \}$$  with the following multiplication law:
$$(z,w)(z^{'},w^{'}) := (z+z',w+e^{z}w')$$
}\end{remark}

\begin{remark} {\rm 
The co-adjoint orbit of $\widetilde\Aff({\bf C})$ in ${\mathfrak g}^*$  passing through $F \in {\mathfrak g}^*$ is denoted by 
$$\Omega_{F} := K(\widetilde {\Aff}({\bf C})) F = \{K(g)F \vert g \in \widetilde \Aff({\bf C})\}$$
Then, (see [D]):
\begin{enumerate}
\item Each point $(\alpha,0,0,\delta)$ is 0-dimensional co-adjoint orbit $\Omega_{(\alpha,0,0,\delta)}$
\item The open set $\beta^{2}+\gamma^{2} \ne $ 0 is the single 4-dimensional co-adjoint orbit $\Omega_{F} = \Omega_{\beta^{2}+\gamma^{2} \ne 0} $. We shall also use $\Omega_{F}$ in form $\Omega_{F} \cong {\bf C} \times {{\bf C}}^*$.
\end{enumerate}
}\end{remark}

\begin{remark}{\rm 
Let us denote:
$$\mathbf H_{k} = \{w=q_{1}+iq_{2} \in {\bf C} \vert -\infty< q_1<+\infty ; 2k\pi < q_2< 2k\pi+2\pi\}; k=0,\pm1,\dots$$
$$L=\{{\rho}e^{i\varphi} \in {\bf C} \vert 0< \rho < +\infty; \varphi = 0\} \mbox{ and } {\bf C} _{k } = {\bf C} \backslash L$$
is a univalent sheet of the Riemann surface of the complex variable multi-valued analytic function $\Ln(w)$, ($k=0,\pm 1,\dots$)
Then there is a natural diffeomorphism $w \in mathbf H_{k} \longmapsto e^{w} \in {\bf C}_k$ with each $k=0,\pm1,\dots.$ Now consider the map:
$${\bf C} \times {\bf C} \longrightarrow \Omega_F = {\bf C} \times {\bf C}^*$$
$$(z,w) \longmapsto (z,e^w),$$
with a fixed $k \in \mathbf Z$. We have a local diffeomorphism 
$$\varphi_k: {\bf C} \times {\bf H}_k \longrightarrow {\bf C} \times {\bf C}_k$$
 $$(z,w) \longmapsto (z,e^w) $$
This diffeomorphism $\varphi_k$ will be needed in the all sequel.
}\end{remark}

On ${\bf C}\times {\bf H}_k$ we have the natural symplectic form 
\begin{equation}\omega = \frac{1 }{2}[dz \wedge dw+d\overline {z} \wedge d\overline {w}],\end{equation} induced from $\mathbf C^2$.
Put $z=p_1+ip_2,w=q_1+iq_2$ and $(x^1,x^2,x^3,x^4)=(p_1,q_1,p_2,q_2) \in {\bf R}^4$, then
$$\omega = dp_1 \wedge dq_1-dp_2 \wedge dq_2.$$ The corresponding symplectic matrix of $\omega$ is 
$$ \wedge = \left(\begin{array}{cccc} 0 & -1 & 0 & 0 \\
		           1 & 0 & 0 & 0 \\
                                                0 & 0 & 0 & 1 \\
                                                0 & 0 & -1 & 0 \end{array}\right)
\mbox{   and   }          
 \wedge^{-1} = \left(\begin{array}{cccc} 0 & 1 & 0& 0 \\
		           -1 & 0 & 0 & 0 \\
                                                0 & 0 & 0 & -1 \\
                                                0 & 0 & 1 & 0 \end{array}\right)$$

We have therefore the Poisson brackets of functions as follows. With each $f,g \in {\bf C}^{\infty}(\Omega)$ 
$$\{f,g\} = \wedge^{ij}\frac{\partial f }{\partial x^i}\frac{\partial g}{\partial x^j} =  
\wedge^{12}\frac{\partial f }{\partial p_1}\frac{\partial g}{\partial q_1}+
\wedge^{21}\frac{\partial f }{\partial q_1}\frac{\partial g}{\partial p_1} +
\wedge^{34}\frac{\partial f}{\partial p_2}\frac{\partial g}{\partial q_2} +
\wedge^{43}\frac{\partial f}{\partial q_2}\frac{\partial g}{\partial p_2} = $$
$$\ \ \ \ \ \ \ =\frac{\partial f }{\partial p_1}\frac{\partial g}{\partial q_1} -
\frac{\partial f}{\partial q_1}\frac{\partial g}{\partial p_1} -
\frac{\partial f}{\partial p_2}\frac{\partial g }{\partial q_2} +
\frac{\partial f}{\partial q_2}\frac{\partial g }{\partial p_2} = $$
$$\ \ =2\Bigl[\frac{\partial f}{\partial z}\frac{\partial g}{\partial w} -
\frac{\partial f}{\partial w}\frac{\partial g}{\partial z} +
\frac{\partial f}{\partial \overline z}\frac{\partial g}{\partial \overline{w}} -
\frac{\partial f}{\partial \overline w}\frac{\partial g}{\partial \overline z}\Bigr]$$  

\begin{proposition}
Fixing the local  diffeomorphism $\varphi_k (k \in {\bf Z})$, we have: 
\begin{enumerate}
\item For any element $A \in \aff(\mathbf C)$, the corresponding Hamiltonian function $\widetilde{A}$  in local coordinates $(z,w)$ of the orbit $\Omega_F$  is of the form
$$\widetilde A\circ\varphi_k(z,w) = \frac{1}{2} [\alpha z +\beta e^w + \overline{\alpha} \overline{z} + \overline{\beta}e^{\overline {w}}]$$
\item In local coordinates $(z,w)$ of the orbit $\Omega_F$, the symplectic Kirillov form $\omega_F$ is just the standard form (1).
\end{enumerate}
\end{proposition}
{\sl Proof}.
 1$^o$  Each element   $ F \in \Omega \subset (\aff({\bf C}))^*$  is of the form $$F = zX^{*} +e^{w}Y^{*} = \left(\begin{matrix}z & 0 \cr e^{w} & 0  \end{matrix}\right)$$ in local Darboux coordinates $(z,w)$. From this implies
$$\widetilde{A}(F) = \langle F,A\rangle = \Re \tr (F.A) = $$
$$=\Re\tr \left(\begin{matrix}\alpha z & \beta z \cr \alpha e ^{w} & \beta e^ {w}\end{matrix}\right) = \frac{1}{ 2} [\alpha z +\beta e^w + \overline{\alpha} \overline{z} + \overline{\beta}e^{\overline {w}}]$$

2$^o$ Using the definition of the Poisson brackets$\{, \}$, associated to a symplectic form $\omega$, we have
\begin{equation}\{\widetilde {A},f \} = \alpha\frac{\partial f }{\partial w} - \beta e^{w} \frac{\partial f} {\partial z} -\overline{\beta} e^ {\overline{w}} \frac{\partial f }{\partial {\overline{z}}}+ \overline{\alpha}\frac{\partial f}{\partial \overline w} \end{equation}
Let us from now on denote by $\xi_A$ the Hamiltonian vector field (symplectic gradient) corresponding to the Hamiltonian function $\tilde{A}$, $A\in \aff(\mathbf C)$. 

Now we consider two vector fields:
$$\xi_A =  \alpha_1\frac{\partial  }{\partial w} - \beta_1 e^{w} \frac{\partial  }{\partial z} -\overline{\beta_1} e^ {\overline{w}} \frac{\partial  }{\partial {\overline{z}}}+ \overline{\alpha_1}\frac{\partial  }{\partial{\overline{w}}} ; \xi_B =  \alpha_2\frac{\partial  }{\partial w} - \beta_2 e^{w} \frac{\partial }{\partial z} -\overline{\beta_2} e^ {\overline{w}}\frac{\partial  }{\partial {\overline{z}}}+ \overline{\alpha_2}\frac{\partial  }{\partial{\overline{w}} },$$
where $A = \left(\begin{matrix}\alpha_1& \beta_1 \cr 0 & 0 \end{matrix}\right); B = \left(\begin{matrix}\alpha_2& \beta_2 \cr 0 & 0 \end{matrix}\right) \in \aff ({\bf C})$.
It is easy to check that
$$\xi_A \otimes \xi_B = \beta_1 \beta_2 e ^{2w} \frac{\partial  }{\partial z} \otimes \frac{\partial }{\partial z}+
 \alpha_1 \alpha_2  \frac{\partial  }{\partial w} \otimes \frac{\partial  }{\partial w}+ 
\overline {\beta_1}  \overline{\beta_2} e ^{2\overline w}\frac{\partial  }{\partial{\overline { z}}} \otimes \frac{\partial  }{\partial{\overline {z}}} + \overline {\alpha_1}  \overline{\alpha_2} \frac{\partial  }{\partial{\overline { w}}} \otimes \frac{\partial  }{\partial{\overline {w}}} + $$
$$+ ( \alpha_1 \beta_2 - \alpha_2\beta_1 )e ^{w}\frac{\partial  }{\partial z} \otimes \frac{\partial }{\partial w}+ ( \overline{\alpha_1} \overline{ \beta_2} -\overline{\alpha_2}\overline{\beta_1} )e ^{\overline{w}}\frac{\partial  }{\partial{\overline  z}} \otimes \frac{\partial  }{\partial{\overline w}}+( \alpha_1 \overline{ \beta_2} -\alpha_2\overline{\beta_1} )e ^{\overline{w}}\frac{\partial  }{\partial{\overline  z}} \otimes \frac{\partial  }{\partial w}+ $$
$$+( \overline{\alpha_1}  \beta_2 -\overline{\alpha_2}\beta_1 )e ^{w}\frac{\partial  }{\partial z} \otimes \frac{\partial  }{\partial{\overline w}}+(\beta_1\overline{\beta_2} - \overline{\beta_1}\beta_2 )e^ {w + \overline{w}} \frac{\partial  }{\partial z} \otimes \frac{\partial  }{\partial{\overline z}}+
(\alpha_1\overline{\alpha_2} - \overline{\alpha_1}\alpha_2 ) \frac{\partial }{\partial w} \otimes \frac{\partial }{\partial{\overline w}}.$$
Thus, $$ \langle\omega,\xi_A \otimes \xi_B\rangle = \frac{1}{2} \Bigl[( \alpha_1 \beta_2 - \alpha_2\beta_1 )e ^{w} + 
( \overline{\alpha_1} \overline{ \beta_2} -\overline{\alpha_2}\overline{\beta_1} )e ^{\overline{w}} \Bigr] = \Re \tr(F.[A,B]) = \langle F,[A,B]\rangle. $$
The proposition  is proved.
\hfill$\square$

\section{Computation of Operators $\hat{\ell}_A^{(k)}$.}

\begin{proposition}\label{Proposition 3.1}
With $A,B \in \aff({\bf C})$, the Moyal $\star$-product satisfies the relation:
\begin{equation} i \widetilde{A} \star i \widetilde{B} - i \widetilde{B} \star i \widetilde{A} = i[\widetilde{A,B} ]\end{equation}
\end{proposition}
{\sl Proof}. 
Consider two arbitrary elements $A=\alpha_1X+\beta_1Y$; $B=\alpha_2X+\beta_2Y$. 
Then the corresponding Hamiltonian functions are:
$$\widetilde{A} = \frac{1}{2} [\alpha_1 z +\beta_1 e^{w} + \overline{\alpha_1}\overline z + \overline{\beta_1} e^ {\overline w}];\widetilde{B} =\frac{1}{2}[\alpha_2 z +\beta_2 e^{w} + \overline{\alpha_2}\overline z + \overline{\beta_2} e^ {\overline w}]$$
It is easy, then, to see that:

$$\begin{array}{rcl}P^0(\widetilde{A},\widetilde{B}) &=& \widetilde{A}. \widetilde{B}\\
P^{1}(\widetilde{A},\widetilde{B}) &=& \{ \widetilde{A}. \widetilde{B} \} = 2\Bigl[\frac{\partial\widetilde{A}}{\partial z}\frac{\partial \widetilde{B}}{\partial w} -
\frac{\partial\widetilde {A}}{\partial w}\frac{\partial\widetilde {B}}{\partial z} +
\frac{\partial\widetilde {A}}{\partial \overline z}\frac{\partial \widetilde {B}}{\partial \overline w} -
\frac{\partial\widetilde {A}}{\partial \overline w}\frac{\partial \widetilde{ B}}{\partial \overline z}\Bigr] \\
 &=& \frac{1}{2} \Bigl[( \alpha_1 \beta_2 - \alpha_2\beta_1 )e ^{w} + 
( \overline{\alpha_1} \overline{ \beta_2} -\overline{\alpha_2}\overline{\beta_1} )e ^{\overline{w}} \Bigr]\\
\mbox{ and  }& & \\
P^{r}(\widetilde{A},\widetilde{B}) &=& 0, \   \forall r \ge 2.\end{array}$$
\par
Thus,
$$ i\widetilde{A} \star i \widetilde{B} - i \widetilde{B} \star i \widetilde{A} =\frac{1 }{2i}\Bigl[ P^1(i\widetilde{A},i\widetilde{B}) - P^1(i\widetilde{B},i\widetilde{A})\Bigr]  = $$ 
$$=  \frac{i }{ 2} \Bigl[( \alpha_1 \beta_2 - \alpha_2\beta_1 )e ^{w} + ( \overline{\alpha_1} \overline{ \beta_2}- \overline{\alpha_2}\overline{\beta_1} )e ^{\overline{w}} \Bigr]$$
on one hand.
\par
On the other hand, because of $[A,B] = (\alpha_1\beta_2-\alpha_2\beta_1)Y$ we have
$$ i[\widetilde{A,B}] = i\langle F,[A,B]\rangle =\frac{i }{2} \Bigl[( \alpha_1 \beta_2 - \alpha_2\beta_1 )e ^{w} + ( \overline{\alpha_1} \overline{ \beta_2} - \overline{\alpha_2}\overline{\beta_1} )e ^{\overline w} \Bigr]$$
The Proposition is hence proved.
\hfill$\square$

For each $A \in \hbox{aff}({\bf C}$), the corresponding Hamiltonian function is 
$$\widetilde{A} = \frac{1}{2} [\alpha  z + \beta e^{w} + \overline{\alpha}\overline z + \overline{\beta} e^ {\overline w}] $$ 
and we can consider the operator  ${\ell}_A^{(k)}$  acting on dense subspace
$L^2({\bf R}^2\times ({\bf R}^2)^*,\frac{dp_1dq_1dp_2dq_2}{(2\pi)^2} )^{\infty}$
 of smooth functions by left $\star$-multiplication by $i \widetilde{A}$, i.e:
${\ell}_A^{(k)} (f) = i \widetilde{A} \star f$. Because of the relation in Proposition 3.1, we have
\begin{corollary}\label{Consequence 3.2}
 \begin{equation}{\ell}_{[A,B]}^{(k)} = {\ell}_A^{(k)} \star {\ell}_B^{(k)} - {\ell}_B^{(k)} \star {\ell}_A^{(k)} := {\Bigl[ {\ell}_A^{(k)},  {\ell}_B^{(k)}\Bigr]}^{\star}\end{equation}
\end{corollary}

From this it is easy to see that, the correspondence $A \in \aff({\bf C}) \longmapsto {\ell}_A^{(k)} = $i$\widetilde {A} \star$. is a representation of the Lie algebra $\aff({\bf C}$) on the space N$\bigl[[\frac{i}{2}]\bigr] $ of formal power series, see [G] for more detail.

Now, let us denote  ${\mathcal F}_z$(f) the partial Fourier transform of the function f from the variable $z=p_1+ip_2$ to the variable $\xi=\xi_1+i\xi_2$, i.e:
$${\mathcal F}_z(f )(\xi,w) = \frac{1 }{2\pi} \iint_{R^2} e^{-iRe(\xi \overline{z})} f(z,w)dp_1dp_2$$

Let us denote by $$\mathcal F_z^{-1}(f )(\xi,w) = \frac{1}{2\pi} \iint_{R^2} e^{iRe(\xi \overline{z})} f( \xi,w)d \xi_1d\xi_2$$ the inverse Fourier transform.

\begin{lemma}\label{Lemma 3.2}     Putting $g = g(z,w) = {\mathcal F}_z^{-1} (f )(z,w)$ we have:
\begin{enumerate}
\item
$$\partial_z g = \frac{i}{2}\overline\xi g \    ;  \partial_z^{r} g = {(\frac{i}{2}\overline\xi)}^r g, r=2,3,\dots $$
\item
 $$ \partial_{\overline z} g = \frac{i}{2}\xi g \    ;  \partial_{\overline z}^{r} g = {(\frac{i}{2}\xi)}^r g, r=2,3,\dots $$
\item 
$${\mathcal F}_z(zg) = 2i\partial_{\overline\xi}{\mathcal F}_z(g) = 2i\partial_{\overline \xi}f \    ;  {\mathcal F}_z(\overline{z}g) = 2i\partial_{\xi}{\mathcal F}_z(g) = 2i\partial_{\xi}f $$    
\item
$$ \partial_w g = \partial_w ({\mathcal F}_z^{-1}(f)) = {{\mathcal F}_z} ^{-1}(\partial_{w}f);\     \partial_{\overline w}g = \partial_{\overline w} ({\mathcal F}_z^{-1}(f)  = {{\mathcal F}_z}^{-1}(\partial_{\overline w}f)$$
\end{enumerate}  
\end{lemma}
{\sl Proof}. 
First we remark that
$$\partial_z = \frac{1 }{2}(\partial_{p_1} - i\partial_{p_2});\quad     \partial_{\overline z} = \frac{1 }{2}(\partial_{p_1} + i\partial_{p_2})$$
we obtain 1$^o$,2$^o$.

3$^o$  $${\mathcal F}_z(zg) = \frac{1}{2\pi} \iint e^{-i(p_1\xi_1 + p_2\xi_2)}p_1g(z,w)dp_1dp_2 +$$
$$+ i\frac{1}{2\pi} \iint  e^{-i(p_1\xi_1 + p_2\xi_2)}p_2g(z,w)dp_1dp_2 $$
$$= i \partial_{\xi_1}{\mathcal F}_z(g) + i^2 \partial_{\xi_2} {\mathcal F}_z(g) = (i \partial_{\xi_1} - \partial_{\xi_2}) {\mathcal F}_z(g) = 2i \partial_{\overline \xi} {\mathcal F}_z(g) = 2i \partial_{\overline \xi}f. $$
$$ {\mathcal F}_z (\overline {z}g)  = \frac{1}{2\pi} \iint e^{-i(p_1 \xi_1 + p_2 \xi_2)} p_1 g(z,w)dp_1 dp_2 - $$
$$-i\frac{1}{2\pi} \iint  e^{-i(p_1 \xi_1 + p_2 \xi_2)}p_2 g(z,w)dp_1 dp_2  = $$
$$ =2i \partial_{\xi} {\mathcal F}_z (g) = 2i \partial_{\xi}f. $$

4$^o$   The proof is straightforward.

The Lemma \ref{Lemma 3.2} is therefore proved.  \hfill$\square$

We also need another Lemma which will be used in the sequel.

\begin{lemma}\label{Lemma 3.3}     With $g = {\mathcal F}_z^{-1}$$(f)($$z,w)$, we have: 
\begin{enumerate}
\item
$$    {\mathcal F}_z(P^0(\widetilde{A},g)) = i(\alpha \partial_{\overline \xi} + \overline {\alpha} \partial_{\xi})f + \frac{1}{ 2} \beta e^w f + \frac{1}{2} \overline {\beta} e^{\overline w} f. $$
\item
$$   {\mathcal F}_z(P^1(\widetilde{A},g)) = \overline {\alpha} \partial_{\overline w}f + \alpha \partial_{w}f - \overline {\beta} e^{\overline w} (\frac{i}{2} \xi)f - \beta e^w (\frac{i}{2}\overline {\xi})f. $$
\item
$$    {\mathcal F}_z(P^r(\widetilde{A},g)) = {(-1)}^r.2^{r-1}[\overline \beta{ e^{\overline w}} (\frac{i}{2}\xi)^r + \beta e^w (\frac{i}{2}\overline \xi)^r]f \ \ \ \ \ \      \forall r \ge 2. $$
\end{enumerate}
\end{lemma}
{\sl Proof}.
Applying Lemma \ref{Lemma 3.2} we obtain
1$^o $ $$     P^0(\widetilde {A},g) = \widetilde {A}.g = \frac{1}{2} [\alpha zg + \beta e^w g + \overline \alpha \overline z g + \overline \beta e^{\overline w} g].$$
Thus,  $$ {\mathcal F}_z(P^0(\widetilde{A},g)) = \frac{1}{2} [\alpha {\mathcal F}_z(zg) + \beta e^w {\mathcal F}_z(g) + \overline \alpha {\mathcal F}_z(\overline z g) + \overline \beta e^{\overline w} {\mathcal F}_z(g)]  =$$ 
$$\frac{1}{2} [2i\alpha \partial_{\overline \xi} {\mathcal F}_z(g) + 2i\overline \alpha \partial_{\xi} {\mathcal F}_z(g) + \beta e^w {\mathcal F}_z(g) + \overline \beta e^{\overline w} \mathcal F_z(g)] = $$
$$=i(\alpha \partial_{\overline \xi} + \overline {\alpha} \partial_{\xi})f + 
\frac{1 }{ 2} \beta e^w f + \frac{1}{2} \overline {\beta} e^{\overline w} f. $$

2$^o$    $$(P^1(\widetilde{A},g)) = \wedge^{12} \partial_{p_1} \widetilde {A} \partial_{q_1}g +  \wedge^{21} \partial_{q_1} \widetilde {A} \partial_{p_1}g +  \wedge^{34} \partial_{p_2} \widetilde {A} \partial_{q_2}g +  \wedge^{43} \partial_{q_2} \widetilde {A} \partial_{p_2}g  $$
$$ = \overline \alpha \partial_{\overline w}g + \alpha \partial_w g - \overline \beta e^{\overline w} \partial_{\overline z}g - \beta e^w \partial_zg. $$
This implies that:
$$ {\mathcal F}_z(P^1(\widetilde{A},g)) = \overline {\alpha} \partial_{\overline w} {\mathcal F}_z(g) + \alpha \partial_{w} {\mathcal F}_z(g) - \overline {\beta} e^{\overline w} \partial_{\overline z} {\mathcal F}_z(g) - \beta e^w \partial_{z} \mathcal F_z(g) =$$ 
$$= \overline {\alpha} \partial_{\overline w}f + \alpha \partial_{w}f - \overline {\beta} e^{\overline w} (\frac{i}{2} \xi)f - \beta e^w (\frac{i}{2}\overline {\xi})f. $$ 

3$^o$     $$(P^2(\widetilde{A},g)) = \wedge^{21} \wedge^{21} \partial_{q_1q_1} \widetilde {A} \partial_{p_1p_1}g + \wedge^{21} \wedge^{43} \partial_{q_1q_2} \widetilde {A} \partial_{p_1p_2}g + \wedge^{43} \wedge^{21} \partial_{q_2q_1} \widetilde {A} \partial_{p_2p_1}g + $$
$$ + \wedge^{43} \wedge^{43} \partial_{q_2q_2} \widetilde {A} \partial_ {p_2p_2}g = \frac{1}{2}\bigl[(\overline \beta e^{\overline w} + \beta e^w - \beta e^w + \overline \beta e^{\overline w} + \overline \beta e^{\overline w} - \beta e^w + \beta e^w + \overline \beta e^{\overline w}) \partial_{\overline z}^2 g + $$
$$ + (\overline \beta e^{\overline w} + \beta e^w + \beta e^w - \overline \beta e^{\overline w} - \overline \beta e^{\overline w} + \beta e^w + \beta e^w + \overline \beta e^{\overline w}) \partial_z^2 g  \bigr]  =$$
$$= 2\overline \beta e^{\overline w} \partial_{\overline z}^2 g + 2\beta e^w \partial_z^2 g.$$
This implies also that:
$$ {\mathcal F}_z(P^2(\widetilde{A},g)) = 2\overline \beta e^{\overline w} {\mathcal F}_z(\partial_{\overline z}^2 g) + 2\beta e^w {\mathcal F}_z( \partial_z^2 g) = 2\overline \beta e^{\overline w} (\frac{i}{2}\xi)^2f + 2 \beta e^w (\frac{i}{2}\overline \xi)^2f. $$
By analogy,   $$P^3(\widetilde {A},g) = (-1)^3[4\overline \beta e^{\overline w} \partial_{\overline z}^3g + 4\beta e^w \partial_z^3g]. $$
$$ {\mathcal F}_z(P^3(\widetilde{A},g)) = (-1)^3.2^2[\overline \beta e^{\overline w} (\frac{i}{2}\xi)^3f + \beta e^w (\frac{i}{2}\overline \xi)^3f]$$ 
and with $r \ge 4$
$$P^r(\widetilde {A},g) = {(-1)}^r.2^{r-1}[\overline \beta{ e^{\overline w}} \partial_{\overline z}^r g + \beta e^w \partial_z^rg]. $$
$${\mathcal F}_z(P^r(\widetilde{A},g)) = {(-1)}^r.2^{r-1}[\overline \beta{ e^{\overline w}} (\frac{i}{2}\xi)^r + \beta e^w (\frac{i}{2}\overline \xi)^r]f.$$

The Lemma \ref{Lemma 3.3} is therefore proved. \hfill$\square$

\begin{proposition}\label{Proposition 3.4}
For each $A = \left(\begin{matrix}\alpha & \beta \cr 0 & 0 \cr\end{matrix}\right) \in \aff({\bf C}) $
 and for each compactly supported $C^{\infty}$-function $f \in C_0^{\infty}({\bf C} \times {\bf H}_k)$, we have:
\begin{equation} {\ell}_A^{(k)}{f} := {\mathcal F}_z \circ \ell_A^{(k)} \circ {\mathcal F}_z^{-1}(f) = [\alpha (\frac{1}{2} \partial_w - \partial_{\overline \xi})f + \overline \alpha (\frac{1 }{2}\partial_{\overline w} - \partial_\xi)f + \end{equation}
$$+\frac{i}{2}(\beta e^{w-\frac{1}{2}\overline \xi} + \overline \beta e^{\overline w - \frac{1}{2} \xi})f] $$
\end{proposition}
{\sl Proof}.
Applying Lemma \ref{Lemma 3.3}, we have:
$${\ell}_A^{(k)}(f):= {\mathcal F}_Z(i \widetilde{A} \star {\mathcal F}_z^{-1} (f))= i{\sum_{r \ge 0}\frac{1}{r!}(\frac{1}{2i})^r \mathcal F_z \Bigl( P ^r( \widetilde{A}, {\mathcal  F}_z^{-1} (f))\Bigr)} = $$
$$=i\Bigl\{[i(  \alpha \partial_{\overline \xi}+ \overline{\alpha}\partial _{\xi}) f + \frac{1}{2} \beta e^{w} f +\frac{1}{2} \overline{\beta} e^{\overline{w}} f ]+\frac{1}{1!} (\frac{1}{2i})[  \overline{\alpha}\partial _{\overline{ w}} f + \alpha\partial _{w}f - \overline{\beta} e^{\overline{w}} (\frac{i}{2}\xi)f -$$
$$- \beta e^{w} (\frac{i}{2} \overline{\xi})f]  + \frac{1}{2!} (\frac{-1}{2i})^22[  \overline{\beta} e^{\overline{w}} (\frac{i}{2}\xi)^2 f + \beta e^{w} ({\frac{i}{2} \overline{\xi})}^2 f]  + \dots + $$
$$+\frac{1}{r!} (\frac{-1}{2i})^r 2^ {r-1}[  \overline{\beta} e^{\overline{w}} ({\frac{i}{2}\xi)}^r f + \beta e^{w} ({\frac{i}{2} \overline{\xi})}^{r} f]+ \dots                         \Bigr\}$$
$$=-(  \alpha \partial_{\overline \xi} - \overline{\alpha}\partial _{\xi}) f+\frac{1}{2}(  \overline{ \alpha} \partial_{\overline w}+ \alpha\partial_{w}) f+ i \Bigl\{ \Bigl[ \frac{1}{2} \beta e^{w} +\frac{1}{2} \overline{\beta} e^{\overline{w}}  -\frac{1}{2}\overline{\beta} e^{\overline{w}} (\frac{1}{2}\xi)- \frac{1}{2}\beta e^{w} (\frac{1}{2} \overline{\xi})\Bigr]f+ $$
$$+  \frac{1}{2}. \frac{1}{2!}\Bigl[  \overline{\beta} e^{\overline{w}} ({\frac{-1}{2}\overline\xi)}^2 + \beta e^{w} ({\frac{-1}{2} \xi)}^2 \Bigr]f+ \dots + 
\frac{1}{2} \frac{1}{k!}\Bigl[\overline{\beta} e^{\overline{w}} ({\frac{-1}{2}\xi)}^{r} + \beta e^{w} ({\frac{-1}{2} \overline{\xi})}^{r} \Bigr]f+\dots\Bigr\} $$ 
$$= \Bigl[\alpha( \frac{1}{2}\partial_{w} - \partial _{\overline{\xi}}) +\overline \alpha( \frac{1}{2}\partial_{\overline w } - \partial_\xi) + \frac{i}{2} \beta e^{w}e^{-\frac{1}{2}\overline{\xi}} + \frac{i}{2} \overline{\beta} e^{\overline {w}}e^{-\frac{1}{2}\xi}\Bigr]f $$
$$= \Bigl[\alpha( \frac{1}{2}\partial_w - \partial _{\overline \xi}) +\overline \alpha( \frac{1}{2}\partial_{\overline w} - \partial_\xi) + \frac{i}{2}({\beta e^{w-\frac{1}{2}\overline \xi}} + \overline{\beta} e^{\overline w-\frac{1}{2}\xi})\Bigr]f $$ 
The Proposition is therefore proved. \hfill$\square$

\begin{remark}\label{Remark 3.5}{\rm Setting new variables  u = $w - \frac{1}{ 2}\overline{\xi}$;$v = w + \frac{1 }{2}{\overline{\xi}}$ we have
\begin{equation}\hat {\ell}_A^{(k)}(f) = \alpha\frac{ \partial f }{\partial u}+ \overline{\alpha}\frac{\partial f }{\partial{\overline{u}}}+ \frac{i }{2}(\beta e^{u}+\overline{\beta}e^{\overline{u}})f \vert_{(u,v)}\end{equation}
i.e $\hat {\ell}_A^{(k)} = \alpha\frac{ \partial }{\partial u}+ \overline{\alpha}\frac{ \partial  }{\partial{\overline{u}}}+ \frac{i }{2}(\beta e^{u}+\overline{\beta}e^{\overline{u}})$,which provides a ( local) representation of the Lie algebra  aff({\bf C}).
}\end{remark}

\section{The Irreducible Representations of $\widetilde{\Aff}({\bf C})$ }
Since $\hat {\ell}_A^{(k)}$ is a representation of the Lie algebra  $\widetilde{\hbox {Aff}} ({\bf C})$, we have:
$$\exp(\hat {\ell}_A^{(k)}) = \exp\bigl(\alpha\frac{ \partial }{\partial {u}}+ \overline{\alpha}\frac{ \partial  }{\partial{\overline{u}}}+ \frac{i }{2}(\beta e^{u}+\overline{\beta}e^{\overline{u}})\bigr)$$ is just the corresponding representation of the corresponding connected and simply connected Lie group $\widetilde\Aff ({\bf C})$.

Let us first recall the well-known list of all the irreducible unitary representations of the group of affine transformation of the complex straight line, see [D] for more details.

\begin{theorem}\label{Theorem 4.1}
Up to unitary equivalence, every irreducible unitary representation of $\widetilde{\hbox {Aff}} ({\bf C})$ is unitarily equivalent to one the following one-to-another nonequivalent irreducible unitary representations:
\begin{enumerate}
\item The unitary characters of the group, i.e the one dimensional unitary representation $U_{\lambda},\lambda \in {\bf C}$, acting in ${\bf C}$ following the formula
$U_{\lambda}(z,w) = e^{{i\Re(z\overline{\lambda})}}, \forall (z,w) \in \widetilde{\Aff} ({\bf C}), \lambda \in {\bf C}.$
\item The infinite dimensional irreducible representations $T_{\theta},\theta \in {\mathbf S}^1$, acting on the Hilbert space $L^{2}(\mathbf R\times \mathbf S ^1)$ following the formula:
\begin{equation}\Bigr[T_{\theta}(z,w)f\Bigl](x) = \exp \Bigr(i(\Re(wx)+2\pi\theta[\frac{\Im(x+z) }{2\pi}])\Bigl)f(x\oplus z),\end{equation}
Where \ $(z,w) \in\widetilde{\Aff}({\bf C})$  ;  $x \in {\bf R}\times {\mathbf S} ^1= {\bf C} \backslash \{0\}; f \in L^{2}({\bf R}\times {\mathbf S} ^1);$
$$ x\oplus z = Re(x+z) +2 \pi i \{\frac{\Im(x+z) }{2\pi}\}$$
\end{enumerate}
\end{theorem}
In this section we will prove the following important Theorem which
is very interesting for us both in theory and practice.
\par
\begin{theorem}\label{Theorem 4.2}
The representation $\exp(\hat {\ell}_A^{(k)})$ of the group $\widetilde{\Aff}({\bf C})$ is the irreducible unitary representation 
$T_\theta$ of $\widetilde{\Aff}({\bf C})$ associated, following the orbit method construction, to the orbit $\Omega$, i.e:
$$\exp(\hat {\ell}_A^{(k)})f(x) = [T_\theta (\exp A)f](x),$$
where $f \in L^{2}({\bf R}\times {\mathbf S} ^1) ; A = \begin{pmatrix}\alpha & \beta \cr 0 & 0 \cr\end{pmatrix} \in \aff({\bf C}) ; \theta \in {\mathbf S}^1 ; k = 0, \pm1,\dots$
\end{theorem} 
{\sl Proof}.
Putting $x = e^u \in {\bf C} \backslash \{0\} =  {\bf R} \times $${\mathbf S}^1 $  and recall that
$$\begin{pmatrix}a & b \cr 0 & 1\cr\end{pmatrix} = \exp(A) = \exp \begin{pmatrix}\alpha & \beta \cr 0 & 0 \end{pmatrix},$$

we can rewrite (7) as following:

$$[T_\theta (\exp A)f](e^u) = \exp \Bigl( i(\Re(\frac{e^\alpha -1}{\alpha}\beta e^u)+2\pi\theta[\frac{\Im e^{u+\alpha}}{2\pi}])\Bigr) f(e^{u \oplus \alpha}),$$
where $$u \oplus \alpha = \Re(u+\alpha)+ 2\pi i\{\frac{\Im(u+\alpha) }{2\pi}\} = u+\alpha - 2\pi i[\frac{\Im(u+\alpha) }{2\pi}]. $$

Therefore, for the one-parameter subgroup $\exp tA$, $t \in {\bf R}$, we have the action formula:

$$\bigl[T_\theta (\exp tA)f\bigr](e^u) = \exp \Bigl(i(\Re{\frac{e^{t\alpha}-1}{\alpha}\beta e^u}+2\pi\theta[\frac{\Im e^{u+t\alpha} }{2\pi}])\Bigr)f(e^{u \oplus t\alpha})$$

By a direct computation:
\begin{equation}\frac{\partial }{\partial t} \bigl([T_\theta (\exp tA)f](e^u)\bigr) = \end{equation}
$$=\frac{\partial }{\partial t} \Bigl(\exp\Bigl(\frac{i}{2}(\frac{e^{t\alpha} -1}{\alpha}\beta e^u + \frac{e^{t\overline \alpha} -1}{ \overline \alpha}\overline \beta e^{\overline u})+ 2\pi\theta i[\frac{\Im{ e^{u+t\alpha}}}{2\pi}]\Bigr)\Bigr) + f(e^{u+t\alpha - 2\pi i [\frac{\Im(u+t\alpha)}{2\pi}]}) $$
$$+\exp \bigl(\frac{i }{2}(\frac{e^{t\alpha} -1}{\alpha}\beta e^u + \frac{e^{t\overline \alpha} -1}{\overline \alpha}\overline \beta e^{\overline u})+  2\pi\theta i[\frac{\Im e^{u+t\alpha}}{2\pi}]\bigr)\frac{\partial }{\partial t}f(e^{u+t\alpha - 2\pi i [\frac{\Im (u+t\alpha)}{2\pi }]}) = $$
$$=\frac{i }{2}(\beta e^{u+t\alpha}+\overline \beta e^{\overline u + t\overline \alpha)} \bigl[T_\theta (\exp tA)f\bigr](e^u) +$$
$$+ \exp \Bigl(i(\Re(\frac{e^{t\alpha} -1}{\alpha}\beta e^u)+2\pi\theta i[\frac{\Im e^{u+t\alpha}}{2\pi}]\Bigr)\alpha e^{u \oplus t\alpha} \frac{\partial f }{\partial u} $$
on one hand.

On the other hand, we have:
\begin{equation}\hat {\ell}_A^{(k)}([T_\theta (\exp tA)f](e^u)=\end{equation}
$$=\alpha\frac{\partial }{\partial u}\bigl([T_\theta (\exp tA)f](e^u)\bigr) + \overline \alpha\frac{\partial }{\partial {\overline u}}\bigl([T_\theta (\exp tA)f](e^u)\bigr) + $$
$$+\frac{i }{2}(\beta e^{u}+\overline \beta e^{\overline u})\bigl[T_\theta (\exp tA)f](e^u)\bigr] = $$
$$=\alpha\frac{i }{2}(\frac{e^{t\alpha}-1}{\alpha}\beta e^u) \exp \bigl( i(\Re(\frac {e^{t\alpha} -1}{\alpha}\beta e^u)+2\pi\theta[\frac{\Im e^{u+t\alpha}}{2\pi}]\bigr)\Bigr)f(e^{u \oplus t\alpha})+ $$
$$+ \alpha \exp \bigl(i(\Re(\frac{e^{t\alpha} -1}{\alpha}\beta e^u)+2\pi\theta[\frac{\Im e^{u+t\alpha} }{2\pi}]\bigr)\Bigr)e^{u \oplus t\alpha}\frac{\partial f }{\partial u} + $$
$$+ \overline \alpha \frac{i }{2}(\frac{e^{t\overline \alpha} - 1}{\overline \alpha}\overline \beta e^{\overline u})\exp \bigl(i(\Re(\frac{e^{t\alpha} -1}{\alpha}\beta e^u)+2\pi\theta[\frac{\Im e^{u+t\alpha} }{2\pi}]\bigr)\Bigr) f(e^{u \oplus t\alpha}) +$$
$$+ \frac{i }{2}(\overline \beta e^{\overline u}+ \beta e^u)[T_\theta(\exp tA)f](e^u) = $$
$$= \frac{i }{2}(\beta e^{u + t\alpha} + \overline \beta e^{\overline u + t\overline \alpha})[T_\theta(\exp tA)f](e^u)+$$
$$+\exp \bigl(i(\Re(\frac{e^{t\alpha} -1}{\alpha}\beta e^u)+2\pi\theta[\frac{\Im e^{u+t\alpha} }{2\pi}]\bigr)\Bigr)\alpha e^{u \oplus \alpha t}\frac{\partial f }{\partial u} $$
From (8) and (9) implies that : 
$$\frac{\partial }{\partial t}[T_\theta(\exp tA)f](x) = \hat {\ell}_A^{(k)}\bigl([T_\theta(\exp tA)f](x)\bigr) \ \ \ \ \  \forall x \in {\bf R} \times {\mathbf S}^1. $$
Remark
$$T_\theta(\exp tA)f](e^u) \vert_{t=0} = \exp(2 \pi i[\frac{\Im e^u }{{2\pi}}]\theta)f(e^{u - 2\pi i[\frac{\Im u }{2\pi}]}) 
= f(e^u). $$
This means that: $\exp(\hat {\ell}_A^{(k)})f(x)$ and $[T_\theta(\exp tA)f](x)$ together are the solution of the Cauchy problem
$$\left\{\begin{array}{rcl}\frac{\partial }{\partial t}u(t,x) &=& \hat {\ell}_A^{(k)}u(t,x);\\ 
u(0,x) &=& id.\end{array}\right. $$
The operator $\hat {\ell}_A^{(k)}$ is behaved good enough, so that the Cauchy problem has an unique solution. From this uniqueness we deduce that
$\exp(\hat {\ell}_A^{(k)})f(x) \equiv [T_\theta(\exp tA)f](x) \             \forall x \in {\bf R} \times {\bf S}^1.$
The Theorem is hence proved. \hfill$\square$

\begin{remark}\label{Remark 4.3} {\rm
We say that a real Lie algebra ${\mathfrak g}$ is in the class $\overline{MD}$ if every K-orbit is of dimension, equal 0 or dim ${\mathfrak g}$. Further more, one proved that
([D, Theorem 4.4]) 
Up to isomorphism, every Lie algebra of class $\overline {MD}$ is one of the following:
\begin{enumerate}
\item Commutative Lie algebras.
\item Lie algebra $\aff({\bf R})$ of affine transformations of the real straight line
\item Lie algebra $\aff({\bf C})$ of affine transformations of the complex straight line.
\end{enumerate}
Thus, by calculation for the group of affine transformations of the real straight line $\Aff({\bf R})$ in [DH] and here for the group affine transformations of the complex straight line $\Aff({\bf C})$ we obtained  a description of the quantum $\overline {MD}$ co-adjoint orbits.
}\end{remark}

\end{document}